\documentclass[12pt,twoside,letter]{article}
\usepackage[dvips]{graphicx}
\usepackage{fancyhdr}
\usepackage{bm}

\usepackage{GoetzArticle}
\usepackage{amssymb}
\usepackage{verbatim}
\usepackage{amssymb}
\usepackage{amsmath}
\usepackage[dvips]{color}
\usepackage{rotating}

\usepackage{rawfonts}
\input prepictex
\input pictexwd
\input postpictex

\newcommand{\commentout}[1]{}
\newcommand{\shortproof}[1]{}

\newcommand{\optionappendix}[2]{#1}

\newcommand{\note}[1]{
        \vspace{5 mm}\par \noindent
        \marginpar{\hspace{-.5cm} \textsc{Note}}
        \framebox{
            \begin{minipage}[c]{0.95 \textwidth}
                \flushleft \tt #1
            \end{minipage}
            }
        \vspace{5 mm}\par}

\newcommand{\Z}{\mathbb{Z}}

\newcommand{\h}{{\cal H}}
\newcommand{\R}{\mathbb{R}}

\newcommand{\N}{\mathbb{N}}

\newcommand{\C}{\mathbb{C}}

\newcommand{\eqa}[1]{\begin{eqnarray}#1 \end{eqnarray}}

\newcommand{\calS}{{\cal S}}

\newcommand{\ga}{ {\mathfrak g} }

\newcommand{\leftright}{\blacktriangleright\hspace{-.15cm}\blacktriangleleft}
\newcommand{\updown}{\begin{array}{c}  \blacktriangledown\\[-.4cm] \blacktriangle
\end{array}}

\def\set#1{\{ #1 \}}
\def\abs#1{| #1 |}
\def\norm#1{\|  #1 \|}

\newcommand{\Rh}{{\widehat\R}}
\newcommand{\RRh}{\R {\times}\Rh}
\newcommand{\RRhd}{\R^d{\times}\Rh^d}

\begin{document}

\title{\bf\vspace{-39pt}
On the invertibility of ``rectangular" bi-infinite matrices and
applications in time--frequency analysis
 }
%
%

\author{\ \\ G\"otz E. Pfander \\ \small School of Engineering and Science,
\\ \small Jacobs University Bremen, \\
\small 28759 Bremen, Germany \\ \small g.pfander@iu-bremen.de\\
}

%
%
\maketitle \thispagestyle{fancy}

%
%


\newcommand{\spacing}[1]{\renewcommand{\baselinestretch}{#1}\large\normalsize}
\spacing{1.1}


\maketitle
\renewcommand{\thefootnote}{\fnsymbol{footnote}}

\renewcommand{\thefootnote}{\arabic{footnote}}
\pagenumbering{arabic} \maketitle


\vspace{1cm}
\section*{ABSTRACT}

\vspace{.2cm}
\begin{center}
  \begin{minipage}{12cm}
Finite dimensional matrices having more columns than rows have no
left inverses while those having more rows than columns have no
right inverses. We give  generalizations of these simple facts to
bi--infinite matrices and use those to obtain density results for
$p$--frames of time--frequency molecules in modulation spaces and
identifiability results for operators with bandlimited
Kohn--Nirenberg symbols.
  \end{minipage}
\end{center}


%
%

\section{Introduction}

Matrices in $\C^{m{\times}n}$ are not invertible if $m\neq n$. To
generalize this basic fact from linear algebra to bi--infinte
matrices, we first associate the quadratic shape of $M\in
\C^{m{\times}n}$, $m=n$, to bi-infinite matrices decaying away from
their diagonals, more precisely, by matrices
$M=(m_{j'j})_{j',j\in\Z^d}$ with $m_{j'j}$ small for $\big|
\|j'\|_\infty - \|j\|_\infty\big|$ large. The rectangular shape of
$M\in \C^{m{\times}n}$, $m<n$, is then taken to correspond to
bi-infinite matrices decaying off wedges which are situated between
two slanted diagonals of slope less than one and which are open to
the left and to the right. In short, for $\lambda
> 1$, we assume $m_{j'j}$ small for $\lambda
\|j'\|_\infty-\|j\|_\infty$ positive and large. To this case, we
associate the symbol $\leftright$. Similarly, $M\in
\C^{m{\times}n}$, $m>n$, corresponds to bi-infinite matrices that
are the adjoints of the $\leftright$ matrices described above. That
is, the case $\updown$ is described by: for $\lambda < 1$, we assume
$m_{j'j}$ small for $-\lambda \|j'\|_\infty+\|j\|_\infty$ positive
and large. In both cases, $\lambda\neq 1$ corresponds to $\frac
n m \neq 1$ in the theory of finite dimensional matrices.

We consider bi--infinite matrices that act on weighted $l^p$ spaces,
$1\leq p\leq \infty$. To illustrate our main result 
we first resort to its simplest
case.

\begin{theorem}\label{theorem:mainresultl-2}
Let $M=(m_{j'j}):l^{2}(\Z)\longrightarrow l^{2}(\Z)$ and
$w:\R^+_0\longrightarrow \R^+_0$ satisfies $w(x)=o(x^{-1-\delta})$,
$\delta > 0$.
\begin{enumerate}
    \item If  $|m_{j'j}|< w(\lambda |j'|-|j|)$ for
    $\lambda |j'|-|j|>0$ and $\lambda>1$, then $M$ has
        no bounded left inverses.
    \item If  $|m_{j'j}|< w(-\lambda |j'|+|j|)$ for
    $-\lambda |j'|+|j|>0$ and $\lambda<1$, then $M$ has
        no bounded right inverses.
\end{enumerate}
\end{theorem}

Note that slanted matrices as covered in \cite{ABK07} and in the
wavelets literature \cite{CDM91,DM93,GMP94,Mic92}, decay off slanted
diagonals, that is, $|m_{j',j}|$ small if $\|\lambda j'-j\|_\infty$
large. Since $\|\lambda j'-j\|_\infty\geq \big|\lambda
\|j'\|_\infty-\|j\|_\infty\big|$, the results in
Section~\ref{section:mainresult} apply in the setting of slanted
matrices as well.

After stating and proving our main result as
Theorem~\ref{theorem:mainresult} in
Section~\ref{section:mainresult}, we illustrate its usefulness in
Section~\ref{section:applications} by applying it in the area of
time--frequency analysis. First, Theorem~\ref{theorem:mainresult} is
used to obtain elementary proofs of density theorems for Banach
frames of Gabor systems and of time--frequency molecules in
so-called modulation spaces \cite{FG97,Fei89}. Second, we discuss
how special cases of Theorem~\ref{theorem:mainresult} have been used
to give necessary conditions on the identifiability of
pseudodifferential operators which are characterized by a
bandlimitation of the operators' Kohn--Nirenberg symbols
\cite{KP06,PW06,PW06b}. The background on time--frequency analysis
that is used throughout Section~\ref{section:applications} is given
in Section~\ref{section:applicationsbackground}.

\section{Non--invertibility of ``rectangular'' bi-infinite
matrices}\label{section:mainresult}

Let $l^p_s(\Z^d)$, $1\leq p \leq \infty$, $s\in\R$, be the weighted
$l^p$-space with norm $ \|\{x_j\}
\|_{l^p_s}=\|\,\{(1+\|j\|_\infty)^s\,x_j\}\, \|_{l^p}$, where $
\|\{x_j\}\|_p=\left(\sum_j |x_j|^p \right)^{\frac 1 p}$ and $
\|\{x_j\}\|_\infty=\sup_j |x_j|$. 

\begin{theorem}\label{theorem:mainresult}
Let $1\leq p_1,p_2,q_1,q_2\leq \infty$,
$\frac 1 {p_1}+\frac 1 {q_1}=1$, $\frac 1 {p_2}+\frac 1 {q_2}=1$,
$r_1,r_2,s_1,s_2\in\R$, and
$M=(m_{j'j}):l_{s_1}^{p_1}(\Z^d)\rightarrow l_{s_2}^{p_2}(\Z^d)$.

\begin{enumerate}
\item  If there exists a $\delta \geq 0$ with
$r_1-s_1+\delta>0$ and $ \frac d {p_2}  + r_1+r_2-s_1+s_2+\delta>0$,
and if there exists $\lambda>1$, $K_0>0$, and a function
$w:\R^+_0\longrightarrow \R^+_0$ with $w(x)=o\left(x^{-(\frac 1
{q_1}+\frac 1 {p_2})d-r_1-r_2+s_1-s_2-\delta}\right)$ and
$$\displaystyle
  |m_{j'j}|\leq w(\lambda \|j'\|_\infty-\|j\|_\infty)\
    (1+\|j\|_\infty)^{r_1} \ (1+\|j'\|_\infty)^{r_2}
    ,\quad \lambda \|j'\|_\infty-\|j\|_\infty > K_0,
$$ then $M$ has no bounded left inverses.

\item  If there exists a $\delta \geq 0$ with  $r_2-s_2+\delta>0$
and $ \frac d {p_1}  + r_1+r_2+s_1-s_2+\delta>0$ and if there exists
$0<\lambda<1$, $K_0>0$ and a function $w:\R^+_0\longrightarrow
\R^+_0$  with $w(x)=o\left(x^{-(\frac 1 {p_1}+\frac 1
{q_2})d-r_1-r_2-s_1+s_2+\delta}\right)$ and
$$\displaystyle
  |m_{j'j}|\leq w(-\lambda \|j'\|_\infty+\|j\|_\infty)\
    (1+\|j\|_\infty)^{r_1} \ (1+\|j'\|_\infty)^{r_2}
    ,\quad -\lambda \|j'\|_\infty+\|j\|_\infty > K_0,
$$ $\lambda,K_0>0$, then $M$ has no bounded right inverses.
\end{enumerate}
\end{theorem}
Clearly, Theorem~\ref{theorem:mainresultl-2} is
Theorem~\ref{theorem:mainresult} for $r_1=r_2=s_1=s_2=0$,
$p_1=q_1=p_2=q_2=2$, and $d=1$. Theorem~\ref{theorem:mainresult} is
a direct consequence of
\begin{lemma}\label{theorem:boundedleftinverse}
Let $1\leq p_1,q_1,p_2\leq \infty$, $\frac 1 {p_1}+\frac 1 {q_1}=1$,
and $M=(m_{j'j}):l^{p_1}(\Z^d)\rightarrow l^{p_2}(\Z^d)$. If there
exists a function $w:\R_0^+\rightarrow\R_0^+$ with
$w(x)=o\left(x^{-(\frac 1 {q_1}+\frac 1
{p_2})d-r_1-r_2-\delta}\right)$ satisfying
$$\displaystyle
  |m_{j'j}|\leq w(\lambda \|j'\|_\infty-\|j\|_\infty)\
    (1+\|j\|_\infty)^{r_1} \ (1+\|j'\|_\infty)^{r_2},\quad \lambda \|j'\|_\infty-\|j\|_\infty >
    K_0\,,
$$
for some constants  $\lambda,K_0,r_1,r_2,\delta$, with $\lambda, K_0
>1$, $\delta \geq 0$, $r_1+\delta>0$, and $ \frac d {p_2}  +
r_1+r_2+\delta>0$, then $M$ has no bounded left inverses.
\end{lemma}

\begin{proof}
We begin with the case $p_1 >1$,  $p_2<\infty$ and show that if
$w:\R_0^+\rightarrow\R_0^+$ satisfies $w(x)=o\left(x^{-(\frac 1
{q_1}+\frac 1 {p_2})d-r_1-r_2-\delta}\right)$, $\delta \geq 0$,
$r_1+\delta>0$ and $ \frac d {p_2}  + r_1+r_2+\delta>0$, then
\eqa{A_{K_1}=
    K_1^{p_2 r_1}\sum_{K\geq K_1}K^{p_2 r_2{+}d{-}1}\left(
    \sum_{k\geq K} k^{d{-}1}\, w(k)^{q_1}\right)^{\frac{p_2}{q_1}}\to 0\,
    \text{ as } K_1\to \infty . \label{equation:DoubleSumW}
}
We set $\widetilde{w}(x)=\sup_{y\leq x} w(y)\in o\left(x^{-(\frac 1
{q_1}+\frac 1 {p_2})d-r_1-r_2-\delta}\right)$ and $v\in C_0(\R^+)$
with $\widetilde{w}(x)\leq v(x)\,x^{-(\frac 1 {q_1}+\frac 1
{p_2})d-r_1-r_2-\delta}$. Then
\eqa{
    \sum_{K\geq K_1+2}K^{p_2 r_2{+}d{-}1}&& \hspace{-1cm}\left(
        \sum_{k\geq K} k^{d-1}\,
        w(k)^{q_1}\right)^{\frac{p_2}{q_1}}
            \leq\sum_{K\geq K_1+1}K^{p_2 r_2{+}d{-}1}
            \left(\sum_{k\geq K+1} k^{d-1}\,
           \widetilde{w}(k)^{q_1}\right)^{\frac{p_2}{q_1}}\notag \\
            &\leq&
                \int_{K_1}^\infty x^{p_2 r_2{+}d{-}1}
                \left(
                \int_x^\infty y^{d-1}\, \widetilde{w}(y)^{q_1}\,dy
                    \right)^{\frac{p_2}{q_1}} \,dx
                \notag \\
 \shortproof{
            &\leq&  \int_{K_1}^\infty x^{p_2 r_2{+}d{-}1}
            \left( \int_x^\infty y^{d-1}\, v(y)^{q_1}
                        y^{-q_1(\frac 1 {q_1}+\frac 1 {p_2})
                             d-q_1 r_2-q_1 r_1-q_1\delta}\,dy
                        \right)^{\frac{p_2}{q_1}} \,dx\notag \\
            }
 \shortproof{
            &\leq&  \int_{K_1}^\infty x^{p_2 r_2{+}d{-}1}
            \left(\int_x^\infty  y^{d-1}\,v(y)^{q_1}
                        y^{-(1+\frac {q_1} {p_2})
                             d-q_1 r_2-q_1 r_1-q_1\delta}\,dy\right)^{\frac{p_2}{q_1}} \,dx\notag \\
            }
            &\leq&  \int_{K_1}^\infty x^{p_2 r_2{+}d{-}1}
            \left(\int_x^\infty  v(y)^{q_1}
                        y^{-1-\frac {q_1} {p_2}
                             d-q_1 r_2-q_1 r_1-q_1\delta}\,dy
                             \right)^{\frac{p_2}{q_1}} \,dx\notag
                             \\
            &\leq&
                    \frac{\|v|_{[K_1,\infty)}\|_\infty^{p_2}}
                        {\frac {q_1} {p_2}
                             d+q_1 r_2+q_1 r_1+q_1\delta}
                        \int_{K_1}^\infty x^{p_2{r_2}{+} d {-} 1}\,
                        x^{-d-p_2 r_2-p_2 r_1-p_2\delta}\,dx
                     \notag\\
 \shortproof{
            &\leq&
                    \frac{\|v|_{[K_1,\infty)}\|_\infty^{p_2}}{\frac {q_1} {p_2}
                             d+q_1 r_2+q_1 r_1+q_1\delta}
                        \int_{K_1}^\infty  x^{-p_2 r_1-p_2\delta-1}\,dx
                  \notag\\
  }
  &\leq&
                \frac{\|v|_{[K_1,\infty)}\|_\infty^{p_2}}{(r_1+\delta)( q_1
                             d+p_2q_1 r_2+p_2q_1 r_1+p_2q_1\delta)}
  K_1^{-p_2 r_1-p_2\delta}
  =o(K_1^{-p_2 r_1}), \notag
} since $\|v|_{[K_1,\infty)}\|_\infty\to 0$ as $K_1\to \infty$ and
(\ref{equation:DoubleSumW}) follows.

To show that $ \inf_{x\in
l_0(\Z^d)}\{\frac{\|Mx\|_{l^{p_2}}}{\|x\|_{l^{p_1}}}\}=0$, we fix
$\epsilon>0$ and note that (\ref{equation:DoubleSumW}) provides us
with a $K_1>K_0$ satisfying $A_{K_1} \leq  (2^d
d)^{-\frac{p_2}{q_1}-1}2^{-p_2 r_2}\left(\frac{\lambda -
1}{\lambda}\right)^{p_2 r_1}\epsilon^{p_2}\,.
$

Set $N=\left\lceil \frac{ \lambda (K_1+1)}{\lambda - 1}
\right\rceil$ and $\widetilde{N}=\lceil \frac{N}{\lambda}
\rceil+K_1.$ Then $\frac{\lambda(K_1+1)}{\lambda - 1}
  \leq N\leq\frac{\lambda(K_1+2)}{\lambda - 1}
$  implies $ \lambda N\geq \lambda K_1 + \lambda +N$ and
$
N\geq
 K_1 + \frac N \lambda +1
    > K_1 + \left\lceil \frac N \lambda\right\rceil=\widetilde N.
 $
Therefore, $(2\widetilde{N}+1)^d<(2N+1)^d$ and the matrix $
\widetilde{M} = ( m_{j'j})_{\|j'\|_\infty\leq
\widetilde{N},\|j\|\leq N}: \C^{(2N+1)^d}\longrightarrow
\C^{(2\widetilde{N}+1)^d}$ has a nontrivial kernel. We now choose
$\widetilde{x}\in \C^{(2N+1)^d}$ with $\|\widetilde{x}\|_{p_1}=1$
and $\widetilde{M}\widetilde{x}=0$ and define $x\in l_0(\Z^{2})$
according to $x_j=\widetilde{x}_j$ if $\|j\|_\infty\leq N$ and
$x_j=0$ otherwise.

By construction, we have $\|x\|_{l^{p_1}}=1$, and $(Mx)_{j'}=0$ for
$\|j'\|_\infty\leq \widetilde{N}$. To estimate $(Mx)_{j'}$ for
$\|j'\|_\infty> \widetilde{N}$, we fix $K > K_1$ and one of the
$2d(2 \big(\lceil\frac{N}{\lambda} \rceil+K\big))^{d-1}$ indices
$j'\in\Z^d$ with $\|j'\|_\infty=\lceil\frac{N}{\lambda} \rceil+K$.
We have $\|\lambda j'\|_\infty \geq N+K \lambda$ and $\lambda
\|j'\|_\infty-\|j\|_\infty \geq K\lambda\geq K$ for all $j\in\Z^d$
with $\|j\|_\infty\leq N$. Therefore
\begin{eqnarray*}
|(Mx)_{j'}|^{q_1}
  &=&    \Big|\sum_{\|j\|_\infty\leq N} m_{j'j}x_j\Big|^{q_1}
  \leq \|x\|^{q_1}_{p_1}
    \sum_{\|j\|_\infty\leq N}\left|m_{j'j}\right|^{q_1}
    \\
  &\leq& (1+\|j'\|_\infty)^{q_1 r_2}\,
  \sum_{\|j\|_\infty\leq N}\,(1+\|j\|_\infty)^{q_1 r_1}w(\lambda \|j'\|_\infty-\|j\|_\infty)^{q_1}\\
  &\leq&(1+\|j'\|_\infty)^{q_1 r_2}\,
    (N{+}1)^{q_1 r_1} \sum_{\|j\|_\infty\geq K}w(\|j\|_\infty)^{q_1}\\
  \shortproof
  {&=&(1+\|j'\|_\infty)^{q_1 r_2}\,
     (N{+}1)^{q_1 r_1}\, \sum_{k\geq K}2d(2k)^{d-1} w(k)^{q_1}\\
     }
  &=&2^dd\,(1+\|j'\|_\infty)^{q_1 r_2}\,
     (N{+}1)^{q_1 r_1}\, \sum_{k\geq K}k^{d-1} w(k)^{q_1}.
\end{eqnarray*}
Finally, we  compute
\eqa{&& \hspace{-1cm} \|Mx\|_{l^{p_2}}^{p_2}
  = \sum_{j'\in\Z^d}|(Mx)_{j'}|^{p_2}
  = \sum_{\|j'\|_\infty \geq \lceil \frac{N}{\lambda} \rceil+K_1}
        |(Mx)_{j'}|^{p_2}\notag \\
  &\leq& (2^d d)^{\frac{p_2}{q_1}}
        \sum_{\|j'\|_\infty \geq \lceil \frac{N}{\lambda} \rceil+K_1}
        (1+\|j'\|_\infty)^{p_2 r_2}(N{+}1)^{p_2 r_1}
        \left(\sum_{k\geq \|j'\|_\infty}
            k^{d-1}\, w(k)^{q_1}\right)^{\frac{p_2}{q_1}}\notag\\
  &\leq& (2^d d)^{\frac{p_2}{q_1}} (N{+}1)^{p_2 r_1}
        \sum_{K\geq \lceil \frac{N}{\lambda} \rceil+K_1}
            2d (2K)^{d-1} (K+1)^{p_2 r_2}
        \left(\sum_{k\geq K}
            k^{d-1}\, w(k)^{q_1}\right)^{\frac{p_2}{q_1}}\notag \\
  \shortproof{
  &\leq& (2^d d)^{\frac{p_2}{q_1}+1}2^{p_2 r_2} (N{+}1)^{p_2 r_1}
        \sum_{K\geq \lceil \frac{N}{\lambda} \rceil+K_1}
         K^{p_2 r_2+d-1} \left( \sum_{k\geq K} k^{d-1}\, w(k)^{q_2}\right)^{\frac{p_2}{q_1}}\notag \\
  }
  &\leq&   (2^d d)^{\frac{p_2}{q_1}+1}2^{p_2 r_2}
        \left(\frac{ \lambda (K_1+2)}{\lambda - 1}+1\right)^{p_2 r_1}
            \sum_{K\geq \lceil \frac{N}{\lambda} \rceil+K_1 }
       K^{p_2 r_2+d-1} \left( \sum_{k\geq K} k^{d-1}\, w(k)^{q_1}\right)^{\frac{p_2}{q_1}}\notag \\
  &\leq&   (2^d d)^{\frac{p_2}{q_1}+1}2^{p_2 r_2}
         \left(\frac{ \lambda}{\lambda - 1}\right)^{p_2 r_1}
        (K_1+3)^{p_2 r_1}
        \sum_{K\geq \lceil \frac{N}{\lambda} \rceil+K_1 }
        K^{p_2 r_2+d-1}
        \left(\sum_{k\geq K} k^{d-1}\,
        w(k)^{q_1}\right)^{\frac{p_2}{q_1}}\notag \\
  &\leq&  \epsilon^{p_2}\, , \notag
}
that is, $\|Mx\|_{l^{p_2}}\leq \epsilon$. Since $\epsilon$ was
chosen arbitrarily and $\|x\|_{l^{p_1}}=1$, we have $ \inf_{x\in
l_0(\Z^2)}\{\frac{\|Mx\|_{l^{p_2}}}{\|x\|_{l^{p_1}}}\}=0$ and $M$ is
not bounded below and has no bounded left inverses.

The cases $p_1=1$ and/or $p_2=\infty$ follow
similarly\optionappendix{.}{ as shown in the appendix.}
\end{proof}

{\it Proof of Theorem~\ref{theorem:mainresult}.}

\noindent {\it Part 1.}\quad Let
$M=(m_{j'j}):l^{p_1}_{s_1}(\Z^d)\rightarrow l^{p_2}_{s_2}(\Z^d)$
satisfy the hypothesis of Theorem~\ref{theorem:mainresult}, {\it
part 1}. Suppose that, nonetheless,
$M=(m_{j'j}):l^{p_1}_{s_1}(\Z^d)\rightarrow l^{p_2}_{s_2}(\Z^d)$ has
a bounded left inverse. This clearly implies that
$$
\widetilde M=(\widetilde m_{j'j})=\big( \ m_{j'j} \,
(1+\|j'\|_\infty)^{s_2}\,
    (1+\|j\|_\infty)^{-s_1}\ \big):l^{p_1}(\Z^d)
    \rightarrow l^{p_2}(\Z^d)
$$
has a bounded left inverse which contradicts
Theorem~\ref{theorem:boundedleftinverse}, since for $\lambda
\|j'\|_\infty-\|j\|_\infty > K_0$, we have
\begin{eqnarray*}
   |\widetilde m_{j'j}|&=&\big|\, m_{j'j}
(1+\|j'\|_\infty)^{s_2}\,
    (1+\|j\|_\infty)^{-s_1} )\,\big|\\
    &\leq&  w(\lambda \|j'\|_\infty-\|j\|_\infty)\
    (1+\|j\|_\infty)^{r_1-s_1}\ (1+\|j'\|_\infty)^{r_2+s_2}
\end{eqnarray*}
with $\delta \geq 0$, $r_1-s_1+\delta>0$, $ \frac d {p_2} +
r_1+{r_2} -s_1+s_2 +\delta>0$, and \\ $w(x)=o\left(x^{-(\frac 1
{q_1}+\frac 1 {p_2})d-r_1-r_2+s_1-s_2-\delta}\right)$.

\vspace{.2cm} \noindent {\it Part 2.}\quad The matrix
$M:l^{p_1}_{s_1}(\Z^d)\rightarrow l^{p_2}_{s_2}(\Z^d)$ has a bounded
right inverse if and only if its adjoint
$M^\ast:l^{p_2}_{s_2}(\Z^d)\rightarrow l^{p_1}_{s_1}(\Z^d)$ has a
bounded left inverse. The conditions on $M$ in
Theorem~\ref{theorem:mainresult}, {\it part 2} are equivalent to the
conditions on $M^\ast$ in Theorem~\ref{theorem:mainresult}, {\it
part 1}. The result follows. \hfill $\square$

\section{Applications}\label{section:applications}

Before stating applications of Theorem~\ref{theorem:mainresult}, we
give a brief account of the concepts from time--frequency analysis
that appear in this section.  For  additional background on
time--frequency analysis and, in particular, Gabor frames, see
\cite{Gro01}.

\subsection{Time--frequency analysis and Gabor frames}
\label{section:applicationsbackground}

The Fourier transform of a function $f\in L^1(\R^d)$,   is given by
$ \widehat f(\gamma)=\int f(x) e^{-2\pi i x\cdot \gamma} \,dx$,
$\gamma\in\Rh^d$, where $\Rh^d$ is the dual group of $\R^d$, and
which, aside of notation, equals $\R^d$. The Fourier transform can
be extended to act unitarily on $L^2(\R^d)$ and isomorphically on
the dual space of Schwarz class functions $\calS(\R^{d})$, that is,
on the space of tempered distributions $\calS'(\R^{d})\supset
\calS(\R^{d})$.

The {\it translation operators} $T_y:\calS(\R^d)\longrightarrow
\calS(\R^d)$, $y\in\R^d$, is given by $(T_{y}f)x=f(x{-}y)$,
$x\in\R^d$, and the {\it modulation operator}
$M_\xi:\calS(\R^d)\longrightarrow \calS(\R^d)$ is given by
$(M_{\xi}f)x=e^{2\pi i x \xi}f(x)$, $x\in\R^d$. Both extend
isomorphically to $\calS'(\R^{d})$, and so do their compositions,
the so-called {\it time--frequency shifts} $\pi(z)=\pi(y,\xi)= T_y
M_\xi$, $z=(y,\xi)\in\RRhd$. Note that the adjoint operator
$\pi(z)^\ast$ of $\pi(z)=\pi(y,\xi)$ is $\pi(z)^\ast=e^{2\pi i y \xi
}\pi(-z)$.


The {\it short--time Fourier transform} $V_g f$ of $f\in L^2(\R^d)
\subseteq \calS'(\R^d)$ with respect to a window function $g\in
L^2(\R^d)\setminus\{0\}$  is
$$
    V_g f(z)=\langle f, \pi(z) g \rangle
        = \int_{\R^d} f(x) \overline{g(x-y)}\, e^{-2\pi i (x{-}y)\cdot\xi}\,
        dx
         ,\quad z=(y,\xi) \in \R^d{\times}\Rh^d.
$$
We have $V_g f\in L^2(\R^d{\times}\Rh^d)$ and $\|V_g f
\|_{L^2}=\|f\|_{L^2}\|g\|_{L^2}$.

A central goal in Gabor analysis is to find $g\in L^2(\R^d) $ and
{\it full rank lattices} $\Lambda=A\Z^{2d}\subset\R^d{\times}\Rh^d$,
$A\in\R^{2d{\times 2d}}$ full rank,  which allow the discretization
of the  formula $\|V_g f \|_{L^2}=\|f\|_{L^2}\|g\|_{L^2}$ in the
following sense: for which $g\in L^2(\R^d) $ and full rank lattices
$\Lambda$ exists $A,B>0$ with
\begin{eqnarray}
A\|f\|^2_{L^2} \leq \sum_{z\in\Lambda} |V_g f(z)|^2 \leq
B\|f\|^2_{L^2},\quad f\in L^2(\R^d)\,.
        \label{equation:discretizedReconstruktion}
\end{eqnarray}
If (\ref{equation:discretizedReconstruktion}) is satisfied, then
$(g,\Lambda)=\{\pi ( z ) g \}_{ z \in\Lambda}$ is called
Gabor frame for the Hilbert space $L^2(\R^d)$. More recently, the
question above has been considered for general sequences $\Gamma$ in
$\RRhd$ in place of full rank lattice $\Lambda$
\cite{BCHL06a,BCHL06b,Gro04}.

%
%
To generalize \eqref{equation:discretizedReconstruktion} to Banach
spaces, we adopt the definition of $p$-frames from \cite{AST01}.

\begin{definition}\label{definition:p-frames}
  The Banach space valued sequence $\{g_j\}_{j\in \Z^d}\subseteq X'$,
  $d\in\N$, is an $l^p_s$--frame for the Banach space $X$,
  $1\leq p \leq \infty$, $s\in\R$, if the analysis operator
  $\displaystyle C_{\mathcal F}: X\longrightarrow l_s^p(\Z^d),
    \quad f\mapsto
        \{\langle f, g_j \rangle\}_{j}
  $
  is bounded and bounded below, that is, if there exists $A,B>0$ with
  \begin{eqnarray}
      A \|f\|_X\leq \| \{\langle f, g_j \rangle\} \|_{l^p_s} \leq B
      \|f\|_X,\quad f\in X\,.
      \label{equation:NormEquivalence}
    \end{eqnarray}
\end{definition}

\vspace{-.5cm} Note that in the Hilbert space case $X=L^2(\R^d)$ and
$l^p_s(\Z^{2d})=l^2(\Z^{2d})$,
\eqref{equation:discretizedReconstruktion} implies that
$\displaystyle C_{\mathcal F}$ has a bounded left inverse, while in
the Banach space case \eqref{equation:NormEquivalence} does not
provide us with a left inverse. Therefore, the existence of a
bounded left inverse for $\displaystyle C_{\mathcal F}$ is included
in the definition of the standard generalization of frames to Banach
spaces
\cite{Chr03,Gro91,FZ98}.

Analogously to Definition~\ref{definition:p-frames}, we include a
generalization of Riesz bases in the Banach space setting.
\begin{definition}
  A sequence
  $\{g_j\}_{j\in \Z^d}\subseteq X$,
  $d\in\N$ is called $l^p_s$--Riesz basis in the Banach space $X$,
  $1\leq p \leq \infty$, $s\in\R$, if the synthesis operator $
    D_{\{g_j\}_j }: l^p_s(\Z^{2d})\longrightarrow X,\quad \{c_j\}_j\mapsto
            \sum_j c_j g_j
$  is bounded and bounded below, that is, if there is $A,B>0$ with
  \begin{eqnarray}
      A \|\{c_j\}_j\|_{l^p_s}\leq \|\sum_j c_j g_j  \|_{X} \leq B
      \|\{c_j\}_j\|_{l^p_s},\quad \{c_j\}_j\in l^p_s(\Z^d).
      \notag 
    \end{eqnarray}
\end{definition}

\vspace{-.5cm} The Banach spaces of interest here are the so--called
modulation spaces \cite{Fei83,FG96,Gro03}. Clearly, $V_g
f(z)=\langle f, \pi(z) g \rangle$, $z \in \R^d{\times}\Rh^d$ is well
defined whenever $g\in \calS(\R^d)$ and $f\in \calS'(\R^d)$ (or vice
versa). This together with $ \|V_g f\|_{L^2}=\|g
\|_{L^2}\|f\|_{L^2}$ in the $L^2$--theory motivates the following.
We let $g=\ga\in\calS(\R^d)$ be an $L^2$--normalized Gaussian, that
is, $\ga(x)=2^{\frac d 4}\,e^{-\pi\|x\|^2_2}$, $x\in\R^d$, and
define the {\it modulation space} $M_s^p(\R^d)$, $ s\in\R,\ 1\leq
p\leq \infty$, by
$$
M_s^p(\R^d)=\{f\in \calS'(\R^d):\ V_{\ga} f\in L_s^p(\RRhd)\}
$$
with Banach space norm
$$
\|f\|_{M_s^p}= \|V_{\ga} f\|_{L_s^p}=\left(\int \big|\,(1+\|z\|)^s\,
V_\ga f(z)\,\big|^p\,dz\right)^{\frac 1 p}<\infty\, ,\quad 1\leq
p<\infty,
$$
and the usual adjustment for $p=\infty$.

\begin{example}\label{example:GaussFrame}
\rm For $\lambda <1$, $(\ga,\,\lambda\Z^{2d})$ is an  $l^2$--frame
for $L^2(\R^d)$ \cite{Lyu92,SW92}. Since $\ga\in\calS(\R^d)\subset
M^1_{t}(\R^d)$ for all $t\geq 0$, Theorem 20 in \cite{Gro04} implies
that in this case $(\ga,\,\lambda\Z^{2d})$ is an $l^p_s$--frames for
$M^p_s(\R^d)$ for $s\in\R$ and $1\leq p \leq \infty$. The Wexler-Raz
identity implies that for $\lambda >1$, $(\ga,\,\lambda\Z^{2d})$ is
an $l^2$--Riesz basis in $L^2(\R^d)$. Hence, $D_{(\ga,\lambda
\Z^{2d})}: \ l^2(\Z^{2d})\longrightarrow L^2(\R^d)$ has a bounded
left inverse of the form $C_{(\widetilde \ga,\lambda \Z^{2d})}$
where the so--called dual function $\widetilde \ga$ of $\ga$
satisfies $\widetilde \ga\in\calS(\R^d)$ \cite{Jan95}. The operator
$C_{(\widetilde \ga,\lambda \Z^{2d})}$ is a bounded operator mapping
$M^p_s(\R^d)$ to $l^p_s(\Z^{2d})$. This implies that
$D_{(\ga,\lambda \Z^{2d})}$ has a left inverse and
$(\ga,\,\lambda\Z^{2d})$ is an $l^p_s$--Riesz basis in $M^p_s(\R^d)$
for $s\in\R$ and $1\leq p \leq \infty$.
\end{example}

\subsection{Density results for Gabor $l^p_s$--frames in modulation spaces}

One of the central results in Gabor analysis is the fact that
$(g,\Lambda)$, $g\in L^2(\R^d)$, cannot be a frame for $L^2(\R^d)$
if the measure of a fundamental domain of the full rank lattice
$\Lambda$ is larger than $1$ \cite{Bag90,Dau90,RS95b}.
Generalizations of this result to general sequences $\Gamma$ in
$\RRhd$ require an alternative definition of density
\cite{BCHL06a,Hei07,Lan67}.

\begin{definition}
Let $Q_R=[-R,R]^{2d}\subseteq \RRhd$ and let $\Gamma$ be a sequence
of points in $\RRhd$.  Then
  \begin{eqnarray}
    D^-(\Gamma)
        =\liminf_{R\to\infty}\inf_{z\in\RRhd}
            \frac{|\Gamma\cap Q_R{+}z|}{(2R)^{2d}} \quad
            \text{and}\quad
            D^+(\Gamma)
        =\limsup_{R\to\infty}\sup_{z\in\RRhd}
            \frac{|\Gamma\cap
            Q_R{+}z|}{(2R)^{2d}}\notag
  \end{eqnarray}
are called {\em lower and upper Beurling density of }$\Gamma$. If
$D^+(\Gamma)=D^-(\Gamma)$, then $\Gamma$ is said to have {\em
uniform density} $D(\Gamma)=D^+(\Gamma)=D^-(\Gamma)$.
\end{definition}

\begin{remark}\rm
The density of a sequence $\Gamma$ does not equal the density of its
range set. For example, the density of the sequence $\{\ldots ,-2,
-2,-1,-1,0,0,1,1,2,2,3,3,\ldots\}$ in $\R$ is 2, while the density
of the range of the sequence, namely of $\Z$, is 1.
\end{remark}

In \cite{CDH99}, it was shown that if $(g,\Gamma)$, $g\in
L^2(\R^d)$, $\Gamma\subseteq \RRhd$, is an $l^2$--frame for
$L^2(\R^d)=M^2_0(\R^d)$, then $1\leq D^-(\Gamma)\leq
D^+(\Gamma)<\infty$, a result that has recently been refined by
Theorem 3 and Theorem 5 in \cite{BCHL06b}. For $l^p_s$--frames for
$M^p_s(\R^d)$, Theorem~\ref{theorem:mainresult} implies

\begin{theorem}\label{theorem:GaborApplication}
Let $1\leq p\leq \infty$, $s\in\R$, and $g\in M^\infty_{2d}$ if $s<0$
and $p\neq\infty$ and $g\in M^\infty_{2d+\delta}$, $\delta>s,0$
else. If $(g,\Gamma)$  is an $l^p_s$--frame for $M^p_s(\R^d)$, then
$D^+(\Gamma)\geq 1$.
\end{theorem}

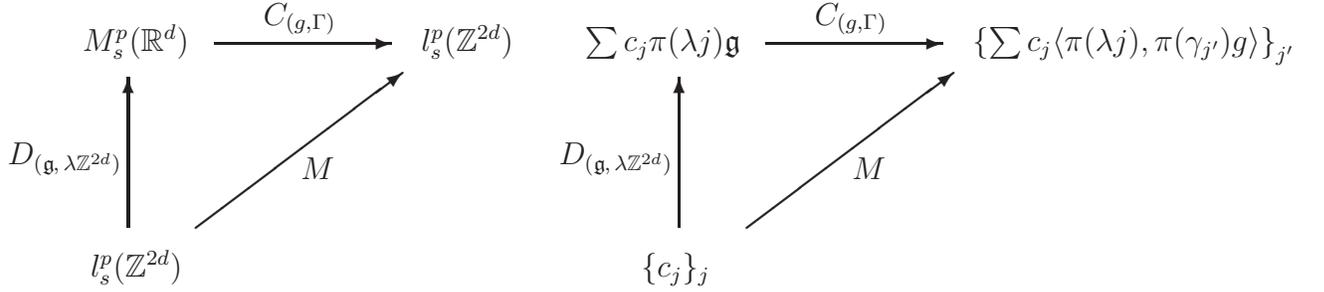
\begin{figure}
\setlength{\unitlength}{.5cm}
\begin{picture}(13,10)\thicklines
    \put(2,7){$M_s^p(\R^d)$}
    \put(6.8,7.7){$C_{(g,\Gamma)}$}
    \put(11,7){$l_s^p(\Z^{2d})$}
    \put(5.5,7.2){\vector(1,0){4.7}}
    \put(2.2,1){$l_s^p(\Z^{2d})$}
    \put(0,4){$D_{(\ga,\,\lambda\Z^{2d})}$}
    \put(3.2,2.3){\vector(0,1){4}}
    \put(5,2.3){\vector(4,3){5.5}}
    \put(7.8,3.6){$M$}
\end{picture}\quad\quad
\begin{picture}(13,10)\thicklines
    \put(.7,7){$\sum c_j \pi(\lambda j) \ga$}
    \put(6.8,7.7){$C_{(g,\Gamma)}$}
    \put(11,7){$\left\{\sum c_j
        \langle \pi(\lambda j),\pi(\gamma_{j'})g\rangle\right\}_{j'}$}
    \put(5.5,7.2){\vector(1,0){4.7}}
    \put(2.2,1){$\set{c_j}_j$}
    \put(0,4){$D_{(\ga,\,\lambda\Z^{2d})}$}
    \put(3.2,2.3){\vector(0,1){4}}
    \put(5,2.3){\vector(4,3){5.5}}
    \put(7.8,3.6){$M$}
\end{picture}

\caption{Sketch of the proof of
Theorem~\ref{theorem:GaborApplication}. We choose $\lambda>1$ so
that $(\ga,\,\lambda \Z^{2d})$ is an $l^p_s$--Riesz basis in
$M^p_s(\R^d)$, so $D_{(\ga,\,\lambda\Z^{2d})}$ is bounded below.
Theorem~\ref{theorem:mainresult} applies to $M=C_{(g,\Gamma)}\circ
D_{(\ga,\,\lambda\Z^{2d})}$, showing that $M$ is not bounded below.
This implies that $C_{(g,\Gamma)}$ is not bounded below and has no bounded left inverses. }
\end{figure}

\begin{proof}
Let $\Gamma$ be given with $D^+(\Gamma)<1$. 
We choose $\lambda>1$ with $1>\lambda^{-4d}> D^+(\Gamma)$ and
$R_0>0$ with
$$|\Gamma\cap
Q_R|< \sup_{z\in\RRhd}
            |\Gamma\cap Q_R{+}z|<\lambda^{-4d}(2R)^{2d},\quad R>R_0.$$
Since $D^+(\Gamma)<\infty$, the sequence $\Gamma$ has no
accumulation points and we can enumerate the sequence $\Gamma$ by
$\Z^{2d}$ so that $\|\gamma_{j'}\|_\infty \leq
\|\gamma_{j''}\|_\infty$ implies $\|{j'}\|_\infty \leq
\|{j''}\|_\infty$ for $j',j''\in\Z^{2d}$. This gives,
$$
\gamma_{j'}\notin Q_R \quad \text{if} \quad
(\,2\|j'\|_\infty-1\,)^{2d}=(\,2(\|j'\|_\infty-1) +1\,)^{2d}\geq
\lambda^{-4d}(2R)^{2d},\quad R>R_0,$$
and, therefore,
\begin{eqnarray}
\gamma_{j'}\notin Q_{\lambda^2\|j'\|_\infty-\frac {\lambda^2}
2}\quad \text{for}\quad \lambda^2\|j'\|_\infty-\tfrac {\lambda^2} 2>
R_0.\label{equation:gammaNotInQ}
\end{eqnarray}

We have
$$
    C_{(g,\Gamma)} \circ D_{(\ga,\,\lambda\Z^{2d})}:l^p_s(\Z^{2d})
        \longrightarrow l^p_s(\Z^{2d}),\quad \{c_j\}_j\mapsto
            \left\{\sum_j c_j \langle \pi(\lambda j) h, \pi(\gamma_{j'})g\rangle
            \right\}=M \{c_j\}_j,
$$
with $M=(m_{j'j})$ and $|m_{j'j}|=|\langle \pi(\lambda j) h,
\pi(\gamma_{j'})g\rangle|=|V_g h(  \gamma_{j'}-\lambda j  )|.$

Note that (\ref{equation:gammaNotInQ}) implies
$$
    \|\gamma_{j'}-\lambda j\|_{\infty}\geq
        \lambda^2\|j'\|_\infty-\tfrac {\lambda^2}
2
        - \|\lambda j\|_\infty
        =\lambda\left(\lambda\|j'\|_\infty - \| j\|_\infty
        -\tfrac {\lambda}
2\right),
$$
and so
$$
|m_{j'j}|=|\langle \pi(\lambda j) \ga,
    \pi(\gamma_{j'})g\rangle|=|V_g \ga(  \gamma_{j'}-\lambda j )|\leq
    w(\lambda\|j'\|_\infty - \| j\|_\infty)
$$
where
$$
    w(\|z\|)=(1+\| z\|)^{-2d-\delta}\sup_{\widetilde z} \big(
(1+\|\widetilde z\|)^{2d+\delta}\ |V_g \ga(\widetilde z)|\big),\quad
z\in\RRhd.
$$
A direct application of Theorem~\ref{theorem:mainresult} implies
that $C_{(g,\Gamma)} \circ D_{(\ga,\,\lambda\Z^{2d})}$ is not
bounded below.  Since $D_{(\ga,\,\lambda\Z^{2d})}$ is bounded below,
we conclude that $C_{(g,\Gamma)}$ is not bounded below which
completes the proof.
\end{proof}

Note that the last lines in the proof of
Theorem~\ref{theorem:GaborApplication} can be modified to apply to
time--frequency molecules which we shall consider in the following.
We say that a sequence
$\{g_{j'}\}_{j'}$ of functions consist of at
$\Gamma=\{\gamma_{j'}\}_{j'}$ $(v,r_1,r_2)$--{\em localized
time--frequency molecules} if
\begin{eqnarray}
     |V_\ga g_{j'}(z)|\leq (1+\|z\|_\infty)^{r_1}
            (1+\|j'\|_\infty)^{r_2} w(\|z-\gamma_{j'}\|_\infty),\quad w=
            o(x^{-v}).\label{equation:molecules}
\end{eqnarray}
If \eqref{equation:molecules} is satisfied for $r_1=r_2=0$, then we
simply speak of at $\Gamma$ $v$--localized time--frequency
molecules. Note that if $\{g_{j'}\}_{j'}\subseteq (M^p_s(\R^d))'$ is
$(v,r_1,r_2)$--localized, then by definition
$\{g_{j'}\}_{j'}\subseteq M^\infty_{v-r_1}(\R^d)$, and,
consequently, if $v-r_1>2d$ we have $\{g_{j'}\}_{j'}\subseteq
M^1(\R^d)$, a fact which we take into consideration when stating the
hypothesis of Theorem~\ref{theorem:MoleculesDplus} and
Theorem~\ref{theorem:MoleculesDminus}

Related concepts of localization were introduced in
\cite{ABK07,Gro04,BCHL06a,BCHL06b}, partly to obtain density results and partly to describe the
time--frequency localization of dual frames of irregular Gabor
frames (see also Remark~\ref{remark:literature}).

\begin{theorem}\label{theorem:MoleculesDplus}
If $\{g_{j'}\}_{j'}\subseteq (M^p_s(\R^d))'\cap M^\infty_{v-r_1}$,
$1\leq p\leq \infty$, $s\in\R$ is an $l_s^p$--frame for
$M^p_s(\R^d)$ which is $(v,r_1,r_2)$--localized at
$\Gamma=\{\gamma_{j'}\}_{j'}$, with $\delta-s,\,
v-r_1-r_2-2d-\delta,\, r_1+\frac{2d}{p}+\delta >0$ and $\delta \geq
0$, then $D^+(\Gamma)\geq 1$.
\end{theorem}

Note that Theorem~9 in \cite{BCHL06b} states that if $\{g_{j'}\}$
is an $l^2$--frame for $L^2(\R^d)$ which consists of at $\Gamma$
$d+\delta$--localized time--frequency molecules, $\delta>0$, then actually
$1\leq D^-(\Gamma)$. Below, we show that components of the proof of
Theorem~\ref{theorem:boundedleftinverse} can be used to obtain some
of the density results given above with $D^+(\Gamma)$ being replaced by
$D^-(\Gamma)$.

\begin{theorem}\label{theorem:MoleculesDminus}
If $\{g_{j'}\}_{j'}\subseteq M^1 (\R^d)$ is an $l^p$--frame for
$M^p(\R^d)$ , $1\leq p\leq \infty$, which is $2d+\delta$--localized
at $\Gamma=\{\gamma_{j'}\}_{j'}$ with $D^+(\Gamma)< \infty$ and
$\delta>0$, then $ D^-(\Gamma)\geq 1$.
\end{theorem}

\begin{proof}
Suppose that $\{g_{j'}\}_{j'}$ is an $l^p_s$--frame for $M^p(\R^d)$
which is  $2d+\delta$--localized at $\Gamma=\{\gamma_{j'}\}_{j'}$,
$D^-(\Gamma)< 1$. For $z_0$, $\alpha_3$ chosen below, we shall
consider the Gabor system $\{\pi(\alpha_3^{-1} j +z_0) \, \ga
\}_{j\in\Z^{2d}}$ which is an $l^p$--Riesz basis for $M^p(\R^d)$. We
shall show that $\{g_{j'}\}$ is not an $l^p$--frame by arguing that
$$\inf_{x\in
    l^p(\Z^d)}\frac{\| C_{\{g_{j'}\}}
        \circ D_{\{\pi(\alpha_3^{-1} j +z_0)\ga\} } x\|_{l^{p}}}
            {\|x\|_{l^{p}}}=0.
$$
To this end, fix $\epsilon>0$. We first assume $1< p < \infty$.

Since $D^+(\Gamma)<\infty$, there exists $\alpha_1\geq 1$ and
$\widetilde {R_0}\geq 1$ with $\infty> \alpha_1^{2d} > D^+(\Gamma)
\geq 0$ and
$$
    |\Gamma\cap Q_R{+}z| \leq  \alpha_1^{2d}
        \, (2R)^{2d}\,,\quad z\in\RRhd,\ R\geq \widetilde {R_0}.
$$
Further, we can pick $\alpha_2, \alpha_3>\tfrac 1 2$ with
$D^-(\Gamma)< \alpha_2^{2d}<\alpha_3^{2d}<1$, and $n_0\in\N$ with
\begin{eqnarray}
  \alpha_2+
        \alpha_1\left( \left(1+\tfrac 1 {n_0}\right)^{2d}-1\right)^{-2d}
    <\alpha_3\left( 1-\tfrac{1}{2n_0} \right)^{2d}
    \notag
    .
\end{eqnarray}

We now choose a monotonically decreasing $w(x)=o(x^{-2d-\delta})$
with $|V_\ga g_{j'}(z)|\leq w(\|z-\gamma_{j'}\|_\infty)$. As
demonstrated in the proof of
Theorem~\ref{theorem:boundedleftinverse}, $w=o(x^{-2d-\delta})$,
$\delta>0$, allows us to pick $\widetilde K_2$ such that for all
$K_2\geq\widetilde K_2$
\begin{eqnarray}
     (2^{2d} 2d)^{\frac{p}{q}+1}
        \sum_{K\geq K_2 }
        K^{2d-1}
        \left(\sum_{k\geq  \frac{\alpha_3}{2\alpha_1} K} k^{2d-1}\,
        w(k)^{q}\right)^{\frac{p}{q}}<\epsilon^p\,. \notag
\end{eqnarray}

Also, there exists $R_0$, $N_0=\lceil\alpha_3 R_0\rceil$,  such that
\begin{itemize}

\item \quad
    $\displaystyle  \text{there exists } z_0\in\RRhd \text{ with }
        |Q_{R_0}{+}z_0\cap \Gamma|\leq \alpha_2^{2d} (2R_0)^{2d}\,;$
    \item \quad
    $\displaystyle R_0\geq \widetilde{R_0}\,n_0$; \quad $\displaystyle N_0 \geq n_0,
    \tfrac{\alpha_1}{\alpha_2}\widetilde{R_0}$;

    \item \quad
    $\displaystyle     (5\tfrac{\alpha_1}{\alpha_3}R_0)^{2d}\,
         w\left(\tfrac {R_0}{n_0}-2\right)<\epsilon\,;$

\item \quad
    $\displaystyle    K_1= N_0 -1-\lceil \alpha_2 N_0 \rceil>1$;
\item \quad
    $\displaystyle    K_2=  2
        \left( \tfrac{\alpha_1}{\alpha_3}N_0
            -\lceil \alpha_2 N_0\rceil\right)
             \geq \widetilde K_2, K_1\,.$
\end{itemize}

The sequence $\Gamma$ has no accumulation point since
$D^+(\Gamma)<\infty$ which implies that we can choose an enumeration
of the sequence $\Gamma$ by $\Z^{2d}$ with $\|{j'}\|_\infty \leq
\|{j''}\|_\infty$ if $\|\gamma_{j'}-z_0\|_\infty \leq
\|\gamma_{j''}-z_0\|_\infty$, $j',j''\in\Z^{2d}$. As mentioned
earlier,  we set $\ga_j=\pi\left( \alpha_3^{-1} j  +z_0 \right)\ga$
for $j\in\Z^{2d}$, and $M=(m_{j'j})=(\langle g_{j'}, \ga_j\rangle)$.

The matrix $ \widetilde{M} = ( m_{j'j})_{\|j'\|_\infty\leq N_0-1,\,
\|j\|\leq N_0}: \C^{(2N_0+1)^d}\rightarrow \C^{(2N_0-1)^d}$ has a
nontrivial kernel, so we may choose $\widetilde{x}\in
\C^{(2N_0+1)^d}$ with $\|\widetilde{x}\|_{p}=1$ and
$\widetilde{M}\widetilde{x}=0$ and define $x\in l_0(\Z^{2})$
according to $x_j=\widetilde{x}_j$ if $\|j\|_\infty\leq N_0$ and
$x_j=0$ otherwise.

To estimate the contributions of  $|(Mx)_{j'}|$ for $j'\in\Z^{2d}$
to $\|Mx\|_{l^p}$,  we consider three cases.

\noindent {\it Case 1.}\quad  $\|j'\|_\infty \leq \lceil\alpha_2
N_0\rceil+ K_1=N_0-1$. \quad This implies $(Mx)_{j'}=0$ by
construction.

\noindent {\it Case 2.}\quad  $\lceil\alpha_2 N_0\rceil + K_1 <
\|j'\|_\infty \leq  \lceil\alpha_2 N_0\rceil+K_2$. \quad Observe
that the set $Q_{R_0+\tfrac {R_0}{n_0}}{+}z_0 \setminus
Q_{R_0}{+}z_0$ consists of a finite number of hypercubes of width
$\tfrac {R_0}{n_0}\geq \widetilde{ R_0}$, so we can estimate
\begin{eqnarray*}
|Q_{R_0+\tfrac {R_0}{n_0}}{+}z_0\cap \Gamma|
    &\leq& \alpha_2^{2d} (2 R_0 )^{2d}
        + \alpha_1^{2d}
        \left(\left( 2\left(R_0+\tfrac {R_0}{n_0}\right)\right)^{2d}
        - \left( 2R_0\right)^{2d}\right) \\
    &\leq& (2R_0)^{2d} \left( \alpha_2^{2d}
        + \alpha_1^{2d} \left(\left(1+\tfrac {1}{n_0}\right)^{2d}
        - 1\right)\right)\\
    &\leq& (2\alpha_3^{-1 } N_0)^{2d} \alpha_3^{2d}
        \left(1-\tfrac 1 {2n_0} \right)^{2d}
        \\
    &\leq&
        \left(2N_0-\tfrac {2N_0} {2n_0} \right)^{2d}
         \leq (2N_0-1)^{2d}
\end{eqnarray*}

Hence, for any $j'$ with $\|j'\|_\infty\geq N_0= \lceil \alpha_2 N_0
\rceil + K_1 + 1$, we have $\gamma_j'\notin Q_{R_0+\tfrac
{R_0}{n_0}}{+}z_0$ and, therefore, for $\|j\|_\infty \leq N_0=\lceil
\alpha_3 R_0 \rceil$ we have
$$
    \| \alpha_3^{-1} j + z_0 -\gamma_{j'} \|_\infty
        =\|(\gamma_{j'}-z_0) - \alpha_3^{-1}j\|_\infty
        \geq R_0+\tfrac
    {R_0}{n_0}-\alpha_3^{-1}\lceil
        \alpha_3 R_0 \rceil\geq \tfrac {R_0}{n_0} - \alpha_3^{-1}
\geq \tfrac {R_0}{n_0} -2,
$$
and, therefore,
\begin{eqnarray}
  |m_{j'j}|=| \langle g_{j'},\ga_j \rangle|
    = |V_\ga g_{j'}(\alpha_3^{-1}j+z_0)|
    \leq  w\left(\| \alpha_3^{-1} j + z_0 -\gamma_{j'} \|_\infty
    \right)
    \leq
     w\left( \tfrac{R_0}{n_0}-2\right)
    .\notag
\end{eqnarray}
This gives
\begin{eqnarray}
&&  \hspace{-3cm}\|Mx|_{ \{j':\,\lceil \alpha_2 N_0\rceil + K_1  <
\|j'\|_\infty
\leq \lceil \alpha_2 N_0\rceil + K_2 \} } \|_p^p \notag \\
    &=&    \sum_{\lceil \alpha_2 N_0\rceil + K_1
                < \|j'\|_\infty \leq \lceil \alpha_2 N_0\rceil + K_2 }
            \big|\sum_{\|j\|_\infty\leq N_0} m_{j'j}x_j\big|^{p} \notag \\
    &\leq&
    \sum_{\lceil \alpha_2 N_0\rceil + K_1
            < \|j'\|_\infty \leq \lceil \alpha_2 N_0\rceil + K_2 }
            \left(\sum_{\|j\|_\infty \leq N_0} |m_{j'j}|^q \right)^{\frac p q}
             \|\widetilde x\|_p^{p}\notag \\
    &\leq&
     w\left(\tfrac
    {R_0}{n_0}-2\right)^p  \hspace{-1cm}
    \sum_{\lceil \alpha_2 N_0\rceil + K_1
            < \|j'\|_\infty \leq \lceil \alpha_2 N_0\rceil + K_2 }
            \hspace{-1cm}
           (2 N_0+1)^{2d\frac p q}\sum_{\|j\|_\infty\leq N_0}  \big|x_j\big|^{p}
           \notag \\
    &\leq&
     w\left(\tfrac
    {R_0}{n_0}-2\right)^p
    (2\cdot 2\tfrac{\alpha_1}{\alpha_3}N_0+1)^{2d}  (2N_0+1)^{2d\frac p q}
    \notag \\ &\leq&
     w\left(\tfrac{R_0}{n_0}- 2 \right)^p
    (5\tfrac{\alpha_1}{\alpha_3}R_0)^{2d(1+\frac p q)}\leq
    \epsilon^p \label{equation:proofDminus1}
\end{eqnarray}

\noindent {\it Case 3.}\quad $\lceil\alpha_2 N_0\rceil+K_2<
\|j'\|_\infty $. \quad For such $j'$, we set  $N=\|j'\|_\infty$ and
obtain $\alpha_1^{-1}(N-\tfrac 1 2)\geq \alpha_1^{-1}(\lceil
\alpha_2 N_0\rceil+K_2+1-\tfrac 1 2)\geq
\tfrac{\alpha_2}{\alpha_1}N_0\geq \widetilde{R_0}$, and, hence,
$$
    |\Gamma \cap Q_{\alpha_1^{-1}(N-\frac 1 2)} +
    z_0 |\, \leq\, \alpha_1^{2d} (2 \alpha_1^{-1}(N-\tfrac 1 2)
    )^{2d}\, =\, (2N-1)^{2d}.
$$
This implies $\gamma_{j'}\notin Q_{\alpha_1^{-1}(\|j'\|_\infty
-\frac 1 2)}+z_0$. Similarly as in {\it Case 2.}, we fix $j'$, $K$
with $\|j'\|_{\infty}= \lceil\alpha_2 N_0\rceil + K$, $K> K_2$, and
conclude that for $\|j\|_\infty \leq N_0$,
\begin{eqnarray*}
   \| \alpha_3^{-1} j + z_0 -\gamma_{j'} \|_\infty
        &=& \|(\gamma_{j'}-z_0) - \alpha_3^{-1}j\|_\infty
        \geq
        \alpha_1^{-1} (\|j'\|_\infty-\tfrac 1 2)-
        \alpha_3^{-1}\|j\|_\infty\\
        &\geq&
        \frac{\alpha_3}{\alpha_1}\|j'\|_\infty -
            \|j\|_\infty - \frac {\alpha_3}{2\alpha_1}
        \\
    &\geq& \frac {\alpha_3}{\alpha_1} \lceil\alpha_2 N_0\rceil
    + 2\frac {\alpha_3}{2\alpha_1}K - N_0 - \frac {\alpha_3}{2\alpha_1}
\\    &\geq& \frac {\alpha_3}{2\alpha_1} \left(K- 2\left(\frac
{\alpha_1} {\alpha_3} N_0 -
        \lceil\alpha_2 N_0\rceil\right)-1\right)+\frac
{\alpha_3}{2\alpha_1}K
    \geq \frac {\alpha_3}{2\alpha_1}K.
\end{eqnarray*}

Therefore,
\begin{eqnarray*}
|(Mx)_{j'}|^{q}
  &=&    \Big|\sum_{\|j\|_\infty\leq N_0} m_{j'j }x_j\Big|^{q}
  \leq \|x\|^{q}_{p}
    \sum_{\|j\|_\infty\leq N_0}\left|m_{j'j}\right|^{q}
  \\&\leq& \,
\sum_{\|j\|_\infty\leq N_0}
    w\left(\frac{\alpha_3}{\alpha_1}\|j'\|_\infty -
            \|j\|_\infty - \frac {\alpha_3}{2\alpha_1}\right)^{q}\\
  &\leq&\,
     \sum_{\|j\|_\infty \geq \frac{\alpha_3}{2\alpha_1}K}
    w(\|j\|_\infty)^{q}
  =
     \sum_{k\geq \frac{\alpha_3}{2\alpha_1}K}2(2d)(2k)^{2d-1} w(k)^{q}\\
  &=&2^{2d}2d\,
      \sum_{k\geq \frac{\alpha_3}{2\alpha_1}K}k^{2d-1} w(k)^{q}.
\end{eqnarray*}
Finally, we compute
\eqa{  \sum_{ \|j'\|_\infty>\lceil\alpha_2 N_0\rceil +K_2 }
\hspace{-1cm}|(Mx)_{j'}|^{p}
  &\leq& (2^{2d} 2d)^{\frac{p}{q}}
        \sum_{\|j'\|_\infty \geq \lceil \alpha_2 N_0 \rceil+K_2}
        \left(\sum_{k\geq \frac{\alpha_3}{2\alpha_1}\|j'\|_\infty}k^{2d-1} w(k)^{q}\right)^{\frac{p}{q}}\notag\\
  &\leq& (2^{2d} 2d)^{\frac{p}{q}}
         \hspace{-.3cm}\sum_{K\geq  \lceil \alpha_2 N_0 \rceil+K_2} \hspace{-.3cm}
            2(2d) (2K)^{2d-1}
        \left(\sum_{k\geq  \frac{\alpha_3}{2\alpha_1} K}
            k^{2d-1}\, w(k)^{q}\right)^{\frac{p}{q}}\notag \\
  &\leq& (2^{2d} 2d)^{\frac{p}{q}+1}
        \hspace{-.4cm}\sum_{K\geq \lceil \alpha_2 N_0 \rceil+K_2} \hspace{-.4cm}
         K^{2d-1} \left( \sum_{k\geq  \frac{\alpha_3}{2\alpha_1} K}
            k^{2d-1}\, w(k)^{q_2}\right)^{\frac{p}{q}}
\leq  \epsilon^p \label{equation:proofDminus2}
}
by hypothesis.  Clearly, \eqref{equation:proofDminus1} and
\eqref{equation:proofDminus2} give $\|Mx\|_{l^{p}}\leq 2^{\frac 1
p}\epsilon$ which completes the proof for $1<p<\infty$. The cases
$p=1$ and $p=\infty$ follow similarly. \end{proof}

\begin{remark}
  \rm If $\{g_j\}=(g,\Gamma)$ and the analysis operator
  $C_{(g,\Gamma)}$ is bounded, then $D^+(\Gamma)<\infty$ follows \cite{CDH99}.
  If
  $\{g_j\}$ are only assumed to be $\Gamma$ localized time--frequency molecules, then
  boundedness of $C_{\{g_j\}}$ does not imply $D^+(\Gamma)<\infty$.
  For example, consider $\{g_j\}=\{\frac 1 {k!}
  \ga\}_{k\in\N}$.
\end{remark}

\begin{remark}\label{remark:literature} \rm
Theorem 9 in \cite{BCHL06b} implies that time--frequency molecules $\{g_j\}$ which are $v$--localized at $\Gamma=\{\gamma_j\}$, $v>d$, and which generate an  $l^2$--frame for $L^2(\R)$ satisfy $1\leq D^-(\Gamma)\leq D^+(\Gamma)$.
Further, Theorem 22 in \cite{Gro04}  states that under the same hypothesis but $v>2d+s$ implies that
being an $l^2$--frame for $L^2(\R^d)$  is equivalent to being an
$l^p_s$--frame for $M^p_s(\R^d)$ for all $1\leq p \leq \infty$ and all $s\geq 0$. This result alone does not imply Theorem~\ref{theorem:MoleculesDplus} nor Theorem~\ref{theorem:MoleculesDminus} as they only assume  that  $\{g_j\}$ is an $l^p_s$--frame for $M^p_s(\R^d)$ for some $p$ and $s$.
Under stronger conditions,  \cite{ABK07} fills this gap.  Namely, Theorem 3.1 and Example 3.1 in \cite{ABK07} show that if    $v>(2d+1)^2+2d$ and $\{g_j\}$ is an  at $\Gamma=\{\gamma_j\}$ $v$--localized  $l^p$--frame for $M^p(\R^d)$ for one $p$, $1\leq p\leq \infty$, then $\{g_j\}$ is  an $l^p$ frame for $M^p(\R^d)$ for all $p$ and therefore for the well studied case $p=2$ \cite{BCHL06b}. This implies  Theorem~\ref{theorem:MoleculesDminus} for $v>(2d+1)^2+2d$.

\end{remark}

\subsection{Identification of operators with bandlimited Kohn--Nirenberg symbols}

A central goal in applied sciences is to identify a partially known
operators $H$ from a single input--output pair $(g,Hg)$. We refer to
an operator class $\h$ as identifiable, if there exists an element
$g$ in the domain of all $H\in\h$ that induces a map
$\Phi_g:\h\longrightarrow Y, \quad H\mapsto Hg$ which is bounded and
bounded below as map between Banach spaces.

In \cite{KP06,PW06}, special cases of
Theorem~\ref{theorem:mainresult} played a crucial role in showing
that classes of pseudodifferential operators with Kohn--Nirenberg
symbol bandlimited to a rectangular domain $[-\tfrac a 2, \tfrac a
2]{\times}[-\tfrac b 2,\tfrac b 2]$ are not identifiable if $ab>1$.
The bandlimitation of a Kohn--Nirenberg symbol to a rectangular
domain $[-\tfrac a 2, \tfrac a 2]{\times}[-\tfrac b 2,\tfrac b 2]$
can be expressed by a corresponding support condition on the
operators so-called {\it spreading function} $\eta_H$\, \footnote{In
fact, the spreading function of an operator is the symplectic
Fourier transform of the operator's Kohn--Nirenberg symbol
\cite{KP06,PW06b}.}. Consequently, we consider operators
$H:D\longrightarrow M^p_s(\R)$, $D\subseteq M^\infty(\R)$, included in
\begin{eqnarray}
\mathcal H^p_s([-\tfrac a 2, \tfrac a 2]{\times}[-\tfrac b 2,\tfrac
b 2])=\left\{ H=\int_{[-\frac a 2, \frac a 2]{\times}[-\frac b
2,\frac b 2]} \eta_H(z)\pi(z)\, dz,\quad \eta_H\in M^p_s(\RRh)
\right\}\, \label{equation:defineoperators}
\end{eqnarray}
and with norm $\|H\|_{\mathcal H^p_s}=\|\eta_H\|_{M^p_s}$. The integral
in \eqref{equation:defineoperators} is defined weakly using $\langle
Hf,h \rangle=\langle \eta_H, V_h f \rangle$\, \footnote{Here,
$\langle \cdot,\cdot\rangle$ is taken to belinear in the first
component and conjugate linear in the second.} \cite{PW06}. In
\cite{KP06} it was shown that

\begin{theorem}
There exists $g\in M^\infty(\R)$ with
   $
        \Phi_g:\h^2_0([-\tfrac a 2, \tfrac a 2]{\times}[-\tfrac b 2,\tfrac
b 2])\longrightarrow M^2_0(\R)$ bounded and bounded below if and
only if $ab\leq 1$.
\end{theorem}

Note that $H^1_0([-\tfrac a 2, \tfrac a 2]{\times}[-\tfrac b
2,\tfrac b 2])$ consists of Hilbert--Schmidt operators,  the norm
$\|\cdot\|_{\mathcal H^2_0}$ is equivalent to the Hilbert--Schmidt
space norm, and $ \| \cdot \|_{M^2_0}$ is a scalar multiple of the
$L^2$--norm.

The main result in \cite{PW06} is
\begin{theorem}
For $ab < 1$ exists $g\in M^\infty(\R)$ with
        $\Phi_g:\h^\infty_0([-\tfrac a 2, \tfrac a 2]{\times}[-\tfrac b 2,\tfrac
b 2])\longrightarrow M^\infty_0(\R)$ bounded and bounded below,
while  for $ab>1$  exists no such $g\in M^\infty(\R)$.
\end{theorem}

Here, we use the generality of Theorem~\ref{theorem:mainresult} to
obtain

\begin{theorem}\label{theorem:identification} Let $1\leq p\leq \infty$
and $s\in\R$. For $ab>1$ exists no $g\in M^\infty(\R)$ with
        $\Phi_g: \mathcal H^p_s
            ([-\tfrac a 2, \tfrac a 2]{\times}[-\tfrac b 2,\tfrac
b 2])\longrightarrow M^p_s(\R)$ bounded and bounded below.
\end{theorem}

\begin{figure} \setlength{\unitlength}{.5cm}
\hspace{-.6cm}\begin{picture}(13,10)\thicklines
    \put(2,7){${\mathcal H}_s^p(\R)$}
    \put(7.7,7.7){$\Phi_g$}
    \put(11,7){$M_s^p(\R)$}
    \put(5.5,7.2){\vector(1,0){4.7}}
    \put(2.2,1){$l_s^p(\Z^{2})$}
    \put(1,4){$D_{\{P_j\}}$}
    \put(3.2,2.3){\vector(0,1){4}}
    \put(5.5,1.3){\vector(1,0){4.7}}
    \put(7.7,1.6){$M$}
    \put(11,1){$l_s^p(\Z^{2})$}
    \put(12.2,6.3){\vector(0,-1){4}}
    \put(12.4,4){$C_{(\ga,\,\lambda\Z^{2d})}$}
\end{picture}\quad \quad
\begin{picture}(13,10)\thicklines
    \put(2,7){$\sum_j c_j P_j$}
    \put(7.7,7.7){$\Phi_g$}
    \put(11,7){$\sum_j c_j P_j g$}
    \put(5.5,7.2){\vector(1,0){4.7}}
    \put(2.2,1){$\left\{ c_j  \right\}_j$}
    \put(1,4){$D_{\{P_j\}}$}
    \put(3.2,2.3){\vector(0,1){4}}
    \put(5.5,1.3){\vector(1,0){4.7}}
    \put(7.7,1.6){$M$}
    \put(11,1){$\left\{\sum c_j \langle P_j g,\pi(\lambda j')\ga\right\}_{j'}$}
    \put(12.2,6.3){\vector(0,-1){4}}
    \put(12.4,4){$C_{(\ga,\,\lambda\Z^{2d})}$}
\end{picture}

\caption{Sketch of the proof of
Theorem~\ref{theorem:identification}. We choose a structured
operator family $\{P_j\}\subseteq \h_s^p$ so that the corresponding
synthesis map $D_{\{P_j\}}:\,\{c_j\}\longrightarrow \sum c_j P_j$
has a bounded left inverse. Further, $C_{(\ga,\,\lambda\Z^{2d})}$
has a bounded left inverse for $\lambda < 1$. We then use
Theorem~\ref{theorem:mainresult} to show that for any $g\in
M^\infty(\R)$, the composition $M = C_{(\ga,\,\lambda\Z^{2d})}\circ
\phi_g \circ D_{\{P_j\}}$ is not bounded below, therefore implying
that $\phi_g: \h_s^p\longrightarrow M^p_s(\R)$ is not bounded below
as well.}
\end{figure}

{\it Sketch of proof.} We assume $a=b$ and $a^2>1$. The general case
$ab>1$ follows similarly.  The goal is to show that for any $g\in
M^\infty(\R)$ which induces a bounded operator $\Phi_g\colon
\h^p_s([-\tfrac a 2, \tfrac a 2]^2) \longrightarrow M^p_s(\R)$, this
operator is not bounded below.

To see this, we pick $\lambda>1$ with $1<\lambda^4<a^2$ and define a
prototype operator $P\in\h^p_s([-\tfrac a 2, \tfrac a 2]^2) $ via
its spreading function  $\eta_P(t,\nu) = \eta(t)\,\eta(\nu)$ where
$\eta$ is smooth,  takes values in $[0,1]$ and satisfies $\eta(t) =
1$ for $\abs{t-a/2}\le a/2\lambda$ and $\eta(t)=
                                       0$ for $\abs{t-a/2}\ge a/2$.

The collection of functions $\{M_{\frac{\lambda}{a} j}\, \eta_P
\}_{j\in\Z^2}$ corresponds to the operator family\\
$\{\pi({\tfrac{\lambda}{a} j})\, P \pi({\tfrac{\lambda}{a} j})^\ast
\}_{j\in\Z^2}$ \cite{PW06}. Further, it forms a Riesz basis for its
closed linear span in $L^2(\RRh)$ and, for $c>0$ sufficiently large,
the collection $\{\pi(\frac{\lambda}{a}j,\frac 1 c k)\, \eta_P
   \}_{j,k\in\Z^2}$
is a frame for $L^2(\R^2)$ \cite{Gro01,Wal92}. Arguing as in
Example~\ref{example:GaussFrame}, we obtain a bounded left inverse
of  $ D_{\{M_{\frac{\lambda}{a} j} \eta_P
\}}:l^p_s(\Z^2)\longrightarrow M^p_s(\RRh)$, thereby showing that
$D_{\{M_{\frac{\lambda}{a} j} \eta_P \}} $ and also the
corresponding operator synthesis map $
D_{\{P_j\}}:l^p_s(\Z^2)\longrightarrow \h^p_s(\RRh)$ with $P_j=
\pi({\tfrac{\lambda}{a} j})\, P \pi({\tfrac{\lambda}{a} j})^\ast$,
$j\in\Z^2$, are bounded below.


For any fixed $g\in M^\infty(\R)$ which induces a bounded map
$\Phi_g:\h^p_s([-\tfrac a 2, \tfrac a 2]^2)\longrightarrow
M^p_s(\R)$ we consider the operator
$$ M=(m_{j j'})=C_{(\ga,\frac{\lambda^2} a)}{\circ}\,\Phi_g{\circ}\, D_{\{P_j\}} \colon
    l^p_s(\Z^2) \longrightarrow l^p_s(\Z^2).$$
We have
$
  \big|m_{j j'}\big|    = \big| \big\langle \pi({\tfrac{\lambda}{a} j}) P
                            \pi({\tfrac{\lambda}{a} j})^\ast\, g,\
                            \pi({\tfrac{\lambda^2}{a} j'})\,\ga
                    \big\rangle\big|
                = \notag \big| V_\ga \, P
                        \pi({\tfrac{\lambda}{a} j})^\ast g
                        \,\big({\tfrac{\lambda}{a} (\lambda j'-j})\big)
                    \big|.
$ In \cite{KP06} it is shown that smoothness and compact support of
$\eta_P$ implies that there exist nonnegative functions $d_1$ and
$d_2$ on $\R$, decaying rapidly at infinity, such that for all $g\in
M^\infty(\R)$, $\abs{Pg(x)}\le \norm{g}_{M^\infty}\,d_1(x)$ and
$\abs{\widehat{Pg}(\xi)}\le \norm{g}_{M^\infty}\,d_2(\xi)$. This
implies that $V_\ga \,P
                        \pi({\tfrac{\lambda}{a} j})^\ast g$ decays rapidly
and independently of $j$, so that we can apply
Theorem~\ref{theorem:mainresult} to show that $M$ is not bounded
below. Since $\tfrac{\lambda^2}{a} < 1$,
Example~\ref{example:GaussFrame} implies that
$C_{(\ga,\frac{\lambda^2} a)}$ is bounded below. Also, $D_{\{P_j\}}$
is bounded below, implying that $\Phi_g$ cannot be bounded below.
Since $ g\in M^\infty(\R)$ was chosen arbitrarily, this completes
the proof. \hfill $\square$

\optionappendix{ }{
\section{Appendix}

\subsection{Proof of Theorem~\ref{theorem:mainresult} for $p_1=1$ and/or $p_2=\infty$}
For $p_1=1$ we proceed analogously as in the case $p_1>1$. We use
$w=o(x^{-\frac d {p_2}-r_1-r_2-\delta})$ to pick $K_1>K_0$ with
$$(K_1+3)^{p_2 r_1}
        \sum_{K\geq \lceil \frac{N}{\lambda} \rceil+K_1 }
        K^{p_2 r_2+d-1}  w(K)^{p_2}
$$
sufficiently small and choose corresponding $N,\widetilde N$. We
estimate
\begin{eqnarray*} |(Mx)_{j'}|
  &=&    \Big|\sum_{\|j\|_\infty\leq N} m_{j'j}x_j\Big|
  \leq \|x\|_{1}
    \sup_{\|j\|_\infty\leq N}\left|m_{j'j}\right|
    \\
  &\leq&\sup_{\|j\|_\infty\leq N}w(\lambda \|j'\|_\infty-\|j\|_\infty) \,
  (1+\|j'\|_\infty)^{r_2}\,(1+\|j\|_\infty)^{r_1}
  \\&\leq& (1+\|j'\|_\infty)^{r_2}\,
  (N{+}1)^{r_1}\sup_{\|j\|_\infty\leq N}w(\lambda \|j'\|_\infty-\|j\|_\infty) \\
  &\leq&(1+\|j'\|_\infty)^{r_2}\,
    (N{+}1)^{r_1}\sup_{\|j\|_\infty\geq K}w(\|j\|_\infty) \\
  &=&(1+\|j'\|_\infty)^{r_2}\,
     (N{+}1)^{r_1}\,w(K),
\end{eqnarray*}
to obtain for $p_2<\infty$,
\eqa{ \|Mx\|_{l^{p_2}}^{p_2}
  &=& \sum_{\|j'\|_\infty \geq \lceil \frac{N}{\lambda} \rceil+K_1}
        |(Mx)_{j'}|^{p_2}\notag \\
  &\leq&
        \sum_{\|j'\|_\infty \geq \lceil \frac{N}{\lambda} \rceil+K_1}
        (1+\|j'\|_\infty)^{p_2 r_2}(N{+}1)^{p_2 r_1}
        w(\|j'\|_\infty)^{p_2}\notag\\
  &\leq&  (N{+}1)^{p_2 r_1}
        \sum_{K\geq \lceil \frac{N}{\lambda} \rceil+K_1}
            2d (2K)^{d-1} (K+1)^{p_2 r_2}
        w(K)^{p_2}\notag \\
  &\leq& 2^{d+p_2 r_2} d (N{+}1)^{p_2 r_1}
        \sum_{K\geq \lceil \frac{N}{\lambda} \rceil+K_1}
         K^{p_2 r_2+d-1} w(K)^{p_2}\notag \\
  &\leq&   2^{d+p_2 r_2} d
         \left(\frac{ \lambda}{\lambda - 1}\right)^{p_2 r_1}
        (K_1+3)^{p_2 r_1}
        \sum_{K\geq \lceil \frac{N}{\lambda} \rceil+K_1 }
        K^{p_2 r_2+d-1}  w(K)^{p_2}<\epsilon^{p_2}.
 \notag
}

For $p_2=\infty$, we use  $w=o(x^{-r_1-r_2-\delta})$ to choose $K_1$
with
$$(K_1+3)^{r_1}
        \sup_{K\geq \lceil \frac{N}{\lambda} \rceil+K_1 }
        K^{r_2}  w(K)\longrightarrow 0 \text{ as } K_1\to\infty$$
        small. Then
\eqa{ \|Mx\|_{l^{\infty}}
  &=& \sup_{\|j'\|_\infty \geq \lceil \frac{N}{\lambda} \rceil+K_1}
        |(Mx)_{j'}|\notag \\
  &\leq&
        \sup_{\|j'\|_\infty \geq \lceil \frac{N}{\lambda} \rceil+K_1}
        (1+\|j'\|_\infty)^{r_2}(N{+}1)^{r_1}
        w(\|j'\|_\infty)\notag\\
  &\leq&  (N{+}1)^{r_1}
        \sup_{K\geq \lceil \frac{N}{\lambda} \rceil+K_1}
             (K+1)^{r_2}
        w(K)\notag \\
  &\leq&
         \left(\frac{ \lambda}{\lambda - 1}\right)^{r_1}
        (K_1+3)^{r_1}
        \sup_{K\geq \lceil \frac{N}{\lambda} \rceil+K_1 }
        K^{r_2}  w(K)<\epsilon.
 \notag
}

 For the remaining case, $p_2=\infty$ and $q_1<\infty$,
$w=o(x^{-\frac d {q_1}-r_1-r_2-\delta})$. We pick $K_1$ with
$$(K_1+3)^{r_1}
        \sup_{K\geq \lceil \frac{N}{\lambda} \rceil+K_1 }
        K^{r_2}  \left(\sum_{k\geq K} k^{d-1}\,
        w(k)^{q_1}\right)^{\frac{1}{q_1}}
        $$
as small as we wish.  With corresponding $N,\widetilde N$, we have
\eqa{ \|Mx\|_{l^{\infty}}
  &=& \sup_{j'\in\Z^d}|(Mx)_{j'}|
  = \sup_{\|j'\|_\infty \geq \lceil \frac{N}{\lambda} \rceil+K_1}
        |(Mx)_{j'}| \notag \\
  &\leq& (2^d d)^{\frac{1}{q_1}}
        \sup_{\|j'\|_\infty \geq \lceil \frac{N}{\lambda} \rceil+K_1}
        (1+\|j'\|_\infty)^{r_2}(N{+}1)^{r_1}
        \left(\sum_{k\geq \|j'\|_\infty}
            k^{d-1}\, w(k)^{q_1}\right)^{\frac{1}{q_1}}\notag\\
  &\leq& (2^d d)^{\frac{1}{q_1}} (N{+}1)^{r_1}
        \sup_{K\geq \lceil \frac{N}{\lambda} \rceil+K_1}
            2d (2K)^{d-1} (K+1)^{r_2}
        \left(\sum_{k\geq K}
            k^{d-1}\, w(k)^{q_1}\right)^{\frac{1}{q_1}}\notag \\
  &\leq& (2^d d)^{\frac{1}{q_1}}2^{r_2} (N{+}1)^{r_1}
        \sup_{K\geq \lceil \frac{N}{\lambda} \rceil+K_1}
         K^{r_2} \left( \sum_{k\geq K} k^{d-1}\, w(k)^{q_2}\right)^{\frac{1}{q_1}}\notag \\
  &\leq&   (2^d d)^{\frac{1}{q_1}}2^{r_2}
        \left(\frac{ \lambda (K_1+2)}{\lambda - 1}+1\right)^{r_1}
            \sup_{K\geq \lceil \frac{N}{\lambda} \rceil+K_1 }
       K^{r_2} \left( \sum_{k\geq K} k^{d-1}\, w(k)^{q_1}\right)^{\frac{1}{q_1}}\notag \\
  &\leq&   (2^d d)^{\frac{1}{q_1}}2^{r_2}
         \left(\frac{ \lambda}{\lambda - 1}\right)^{r_1}
        (K_1+3)^{r_1}
        \sup_{K\geq \lceil \frac{N}{\lambda} \rceil+K_1 }
        K^{r_2}
        \left(\sum_{k\geq K} k^{d-1}\,
        w(k)^{q_1}\right)^{\frac{1}{q_1}}<\epsilon \notag.
}

}


\small\bibliography{../Bibliography/gabor_goetz}

\bibliographystyle{elsart-num}


\end{document}